\documentclass[reqno]{amsart}
\usepackage{amsmath}
\usepackage{amsthm}
\usepackage{amssymb}
\usepackage{epsfig}
\usepackage{stmaryrd}

\usepackage{epic, eepic}

\newcommand{\NN}{\mbox{$\mathbb{N}$}}
\newcommand{\p}{\mbox{$\mathcal{P}$}}

\newcommand{\s}{\mbox{$\sigma$}}
\renewcommand{\t}{\mbox{$\tau$}}
\newcommand{\lsr}{\mbox{$\langle \sigma \rangle$}}
\newcommand\lsrt{\ensuremath{[\lsr,\t]}} 
\newcommand{\lpr}{\mbox{$\langle \pi \rangle$}}

\newcommand\ol{\overline}
\newcommand\sym{\mathcal{S}}

\newcommand{\ifset}[2]{\mbox{\ensuremath{\llbracket {#1}, {#2} \rrbracket}}} 
\newcommand\lsrtt{\ifset{\lsr}{\t}}

\newcommand\la{\langle}
\newcommand\ra{\rangle}

\newcommand\C{\ensuremath{C}}
\def\rp{\ensuremath{\mathrm{RP}}}
\DeclareMathOperator{\rk}{\mathrm{rk}}

\theoremstyle{plain}
\newtheorem{thm}{Theorem}[section]
\newtheorem{cor}[thm]{Corollary}
\newtheorem{conj}[thm]{Conjecture}
\newtheorem{lem}[thm]{Lemma}

\theoremstyle{definition}
\newtheorem{defn}[thm]{Definition}
\newtheorem{ex}[thm]{Example}
\theoremstyle{remark}
\newtheorem{rem}[thm]{Remark}

\title{The M\"obius Function of the Permutation Pattern Poset}

\author{Einar Steingr\'imsson} \address{The Mathematics Institute,
  School of Computer Science, Reykjav\'ik University, Iceland}
\email{einarst@ru.is}

\author{Bridget Eileen Tenner}
\address{Department of Mathematical Sciences, DePaul University, Chicago, Illinois, USA}
\email{bridget@math.depaul.edu}

\thanks{The first author was supported by grants no.\ 060005013 and
  090038011 from the Icelandic Research Fund.}



\begin{document}

\begin{abstract}
  A permutation $\t$ contains another permutation $\s$ as a pattern if
  $\t$ has a subsequence whose elements are in the same order with
  respect to size as the elements in $\s$.  This defines a partial
  order on the set of all permutations, and gives a graded poset $\p$.
  We give a large class of pairs of permutations whose intervals
  in~$\p$ have M\"obius function $0$.  Also, we give a solution to the
  problem when $\s$ occurs precisely once in $\t$, and $\s$ and $\t$
  satisfy certain further conditions, in which case the M\"obius
  function is shown to be either $-1$, $0$ or $1$.  We conjecture that
  for intervals $[\s,\t]$ consisting of permutations avoiding the
  pattern 132, the magnitude of the M\"obius function is bounded by
  the number of occurrences of $\s$ in $\t$.  We also conjecture that the M\"obius function of the
  interval $[1,\t]$ is $-1$, $0$ or $1$.
\end{abstract}

\maketitle
\thispagestyle{empty}

\section{Introduction}

In this paper, permutations are of the letters ${1,2,\dots,n}$ (for
various $n$) and a \emph{pattern} in a permutation $\pi$ is a
subsequence in $\pi$, the relative sizes of whose elements come in
some prescribed order.  For example, a $123$-pattern is simply an
increasing subsequence of length three.  The permutation $246153$ has
two \emph{occurrences} of the pattern $123$, namely the subsequences
$\la246\ra$ and $\la245\ra$.  The permutation $246153$ also has one
occurrence of the pattern $3142$, namely the subsequence $\la4153\ra$.
It is easy to see that the set of all permutations (of arbitrary,
positive, numbers of letters) forms a graded poset $\p$ with respect
to pattern containment.  That is, a permutation $\s$ is smaller than
another permutation $\t$ in $\p$ if $\s$ occurs as a pattern in $\t$.

A classical problem for any poset, first explicitly mentioned by Wilf
\cite{wilf} in this particular case, is to understand its M\"obius
function $\mu$.  The first result so far in this
direction was given by Sagan and Vatter \cite{sagan-vatter}, who
solved the problem in the case of \emph{layered} permutations.  They
noted that the poset of layered permutations is isomorphic to a
certain poset of compositions of an integer, and they gave a formula
for computing the M\"obius function of that poset.  A permutation is
layered if it is the concatenation of decreasing sequences, where the
letters in each sequence are smaller than all letters in later
sequences.  An example of such a permutation is $321546987$.

Looking at examples of the M\"obius function of intervals in the
pattern poset, it is easy to infer that this is hard to understand in
the general case. In particular, the M\"obius function for this poset
does not alternate in sign with rank, which is a discouraging starting
point.

In this paper, we present some results for this problem.  Most of our
results stem from looking at particular occurrences of a pattern $\s$
in a permutation $\t$ and analyzing the complement in $\t$ of an
occurrence.  In particular, in Section~\ref{section:interval blocks},
we give a large class of pairs of permutations $\s$ and $\t$ for which
the M\"obius function of the interval $[\s,\t]$ is 0
(Theorem~\ref{thm:interval blocks}).  These intervals occur when there
is a sequence of consecutive letters of $\t$ that are not in any
occurrence of $\s$, and which form an interval of values
$\{a,a+1,\dots,a+b\}$ for some $a$ and some $b\ge1$.  One such pair is
$(132,859324617)$: no letter of the contiguous subsequence $324$
belongs to any occurrence of the pattern $132$ in $859324617$.

In Section~\ref{section:fixed occurrence}, we solve the problem in the
case where $\s$ occurs precisely once in $\t$ and the complement of
$\s$ in $\t$ satisfies certain conditions.  The M\"obius function in
this case is either $-1$, $0$, or~$1$ (Corollary~\ref{cor:single
  occurrence}).  This follows from a more general result
(Theorem~\ref{thm:fixed occurrence}) about a particular kind of fixed
occurrence of $\s$ in $\t$.  We believe that there are other
interesting classes of intervals also having M\"obius function $-1$,
$0$ or $1$, as we explain below.

In the final section, we conjecture that if the permutations in an
interval $[\s,\t]$ avoid the pattern $132$ (equivalently, if $\t$
avoids $132$), then the absolute value of the M\"obius function of
$[\s,\t]$ is bounded by the number of occurrences of $\s$ in $\t$.
Due to the symmetries among $132$, $231$, $213$ and $312$, the pattern
$132$ can be replaced by any one of the other three patterns in this
conjecture. Thus, this conjecture is equivalent to the conjecture that
if the absolute value of the M\"obius function of an interval
$[\s,\t]$ exceeds the number of occurrences of $\s$ in $\t$ then $\t$
must contain all of the patterns $132$, $231$, $213$ and $312$.
Observe also that if this conjecture holds, then, for intervals
$[\s,\t]$ where $\s$ occurs precisely once in $\t$, the M\"obius
function must be $0$, $1$ or $-1$.

Finally, we conjecture that if a permutation $\t$ avoids $132$, then
the interval $[1,\t]$ has M\"obius function $-1$, $0$ or $1$.

\section{Definitions and notation}

Let $\sym_n$ be the set of permutations of the letters
$\{1,2,\ldots,n\}$.  We represent permutations in one-line notation,
meaning that the permutation $\s \in \sym_k$ is denoted $\s =
\s(1)\s(2) \cdots \s(k)$.  Throughout this section, fix permutations
$\s \in \sym_k$ and $\t \in \sym_n$, with $k \le n$.

\begin{defn}\label{defn:std-form}
If $\pi$ is a permutation of a set of $k$ integers, then its
\emph{standard form} is the permutation of $\{1,2,\ldots,k\}$ whose
letters are in the same relative order of size as those of $\pi$.  We
say that two permutations are \emph{order isomorphic}
if they have the same standard form.
\end{defn}

For example, the standard form of both $3615$ and $4725$ is $2413$,
meaning that $3615$, $4725$, and $2413$ are all order isomorphic to each other.

\begin{defn}\label{defn:pattern}
The permutation $\t$ \emph{contains $\s$} 
(also, $\t$ \emph{has a $\s$-pattern}) if there exist indices $i_1 <
\cdots < i_k$ such that $\s$ is order isomorphic to $\t(i_1) \cdots
\t(i_k)$.  If $\t$ does not contain $\s$, then $\t$ \emph{avoids}~$\s$.
\end{defn}

If $\t$ has a $\s$-pattern, with $i_1 < \cdots < i_k$ as in
Definition~\ref{defn:pattern}, then $\t(i_1) \cdots \t(i_k)$ is an
\emph{occurrence} of $\s$ in $\t$.  The substring $\t(i_{j_1}) \cdots
\t(i_{j_\ell})$ will be denoted $\la \s(j_1) \cdots \s(j_\ell)
\ra$, and occurrences of $\s$ will be distinguished by subscripts:
$\lsr_i$.  
When speaking of an occurrence $\lsr$ it is necessary to know in what larger permutation this occurrence sits ($\t$ in the current discussion).  In many instances this will be clear from the context, and will not be specified further.

\begin{ex}\label{ex:occurrence}
  Let $\t = 74136825$ and $\s = 1243$.  Then $\la1365\ra$ and $\la1385\ra$ are the
  only two occurrences of $\s$ in $\t$.  We can name $\lsr_1 = \la1365\ra$
  and $\lsr_2 = \la1385\ra$.
\end{ex}

If $\t$ has a $\s$-pattern, and $\lsr$ is a particular occurrence of $\s$ in $\t$ with $x \in \t \setminus \lsr$ a letter of $\t$ not in $\lsr$, then we will simplify notation somewhat and write
\begin{equation}\label{eqn:lsr+x}
  \lsr+x := \la \lsr \cup \{x\}\ra
\end{equation}
to indicate the occurrence of a pattern in $\t$ (different from the
$\s$-pattern) formed by the letters $\lsr \cup \{x\}$.  The notation
$\lsr + S$ is analogously defined, whenever $S$ is a subset of letters
in $\t \setminus \lsr$.

\begin{ex}
With the notation of Example~\ref{ex:occurrence}, $\lsr_1 + 2 = \la13625\ra$ and $\lsr_1+7 = \la71365\ra$.
\end{ex}

\begin{defn}
If $\t$ contains $\s$, then write $\s\le\t$. If $\s\le\t$, but $\s$
and $\t$ are not order isomorphic (in particular, $\s\ne\t$), then
write $\s<\t$.  If $\s<\t$ and $k = n-1$, then $\t$ \emph{covers}
$\s$.  Let $\p$ denote the poset of all permutations of arbitrary,
positive, numbers of letters, ordered by pattern containment.  That
is, the partial ordering on $\p$ is defined by the relation $\le$.
\end{defn}

Two elements $s$ and $t$ in a poset $P$ determine an
\emph{interval} $[s,t] = \{x \in P \;|\;s\le x \le
t\}$, with $[s,t] = \emptyset$ if $s \not\le t$ in $P$.  We also define the
half-open interval $[s,t) = \{x \in P \;|\;s \le x < t\}$.

It is easy to see that the permutation pattern poset $\p$ is
\emph{graded}.  That is, the lengths of all maximal chains in an
interval are the same.  More precisely, the length of a maximal chain
between two permutations $\s$ and $\t$, where $\s\le\t$, is the
difference between the number of letters of $\t$ and of $\s$.  This
number is also called the \emph{rank} of the interval.

\begin{defn}\label{defn:mu defn}
  Given an interval $[s, t]$ in a poset $P$, the \emph{M\"obius
    function} $\mu_P$ of this interval, written $\mu$ when no
  confusion will arise, is recursively defined by
\begin{equation}\label{eqn:mu defn}
\mu(s,t) = \begin{cases}
0 & \text{if } s \not\le t,\\
1 & \text{if } s = t, \text{ and}\\
- \sum\limits_{s \le x < t} \mu(s, x) & \text{otherwise.}
\end{cases}
\end{equation}
\end{defn}

\begin{figure}[htbp]
\scalebox{.4}{\input{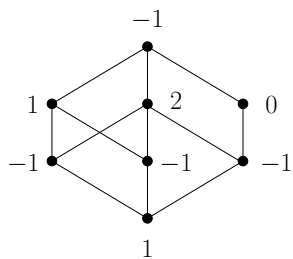}}
\caption{\label{fig:mobius comp} Computing the M\"obius function of an
  interval: the number at each element in this poset gives the value
  of the M\"obius function of the interval from the minimum element to
  that element.}
\end{figure}

From Definition \ref{defn:mu defn}, it is clear that if $s \neq t$,
then the sum of the values $\mu(s, x)$ over all $x$ in the interval
$[s, t]$ equals zero.  For an example of how the M\"obius function can
be computed from Definition \ref{defn:mu defn}, see Figure
\ref{fig:mobius comp}.

\begin{defn}
We regard patterns (permutations) $\s$ in the permutation pattern
poset $\p$ as functions from $\p$ to $\NN$, where $\s(\t)$ is the
number of occurrences of $\s$ in $\t$.
\end{defn}

Note that if $\s\not\le\t$, then $\s(\t)=0$.

\begin{ex}\label{ex:number of patterns}
  Let $\s = 231$ and $\t = 23541$.  Then $\s(\t) = 5$ because there
  are five distinct occurrences of the pattern $\s$ in $\t$: $\la231\ra$,
  $\la251\ra$, $\la241\ra$, $\la351\ra$, and $\la341\ra$.
\end{ex}

\begin{defn}
  If a permutation $\t$ maps a non-singleton interval onto an
  interval, then this image is an \emph{interval block}.  More
  precisely, fix an integer $b \ge 1$.  If $\t$ has a factor
  (consecutive substring) $I = \t(a) \t(a+1) \cdots \t(a+b)$, where
  the set of values $\{\t(a), \t(a+1), \ldots, \t(a+b)\}$ consists of
  all the numbers $a',a'+1,\ldots, a'+b$ for some $a'$, then $I$ is an
  \emph{interval block}.
\end{defn}

\begin{ex}
  The permutation $71342865$ has interval blocks 34, 342, 1342, 65 and 71342865.
\end{ex}

This definition may bring to mind \emph{simple}
permutations: $\pi$ is simple if it has no interval
blocks other than $\pi$ itself. For information about these
permutations, see~\cite{brignall} and~\cite{albert-atkinson}.

\begin{defn}
  Suppose that $\s \le \t$ in $\p$, and fix an occurrence $\lsr$ of
  $\s$ in $\t$.  If $\t$ has an interval block that does not intersect the occurrence
  $\lsr$, then the pair $(\lsr, \t)$ has an \emph{interval block}.
  Otherwise the pair is \emph{interval free}.
\end{defn}

\begin{defn}
Suppose that $\s \le \t$.  If $\t$ has an interval block that does
not intersect any occurrence of $\s$, then the pair $(\s,\t)$ has an
\emph{interval block}.
\end{defn}

\begin{ex}
Let $\s = 2341$ and $\t = 162395784$.  The pair $(\s, \t)$ has
interval block $I = 23$ because the occurrences of $\s$ in $\t$,
namely $\la6784\ra$ and $\la5784\ra$, are each disjoint from $I$.
\end{ex}

To compute the M\"obius function $\mu(\s,\t)$, one examines the
interval $[\s,\t] \subset \p$.  For our purposes, it
will be helpful first to examine a slightly different poset, described
below.  

\begin{defn}\label{defn:occurrence interval}
  Suppose that $\s \le \t$, and that $\lsr$ is a particular occurrence
  of $\s$ in $\t$.  Set $\C := \t \setminus \lsr$ to be the complement of $\lsr$ in $\t$.  
  Let the \emph{occurrence poset} $[\lsr, \t]$ denote the partially ordered set
  consisting of permutations formed by deleting arbitrary subsets of
  letters of $\C$ from $\t$.
  Two such permutations are considered equivalent if they are order isomorphic and if the specified occurrence $\lsr$ is in the same positions in each.  A permutation $\rho$ is covered by $\pi$ if $\rho$ can be obtained from $\pi$ by removing a letter of $\pi \setminus \lsr$.    This defines the partial ordering in $\lsrt$.
\end{defn}

In particular, if $\s(\t) = 1$, then $[\lsr, \t] = [\s,\t]$.  If $\t=\ol63\ol45\ol21$, where we have overlined the letters
in a particular occurrence $\lsr$ of the pattern $321$, then the
permutations obtained by deleting any single letter of $\C=\{3,5,1\}$
are $\ol5\ol34\ol21$, $\ol53\ol4\ol21$ and $\ol52\ol34\ol1$.  Note
that the first two of these are identical as permutations, but are
considered distinct here because $\lsr$ is in
different positions.  Figure \ref{fig:posets} depicts the intervals
$[\s,\t]$ and $[\lsr,\t]$ for this example.

\begin{figure}[htbp]
\scalebox{.4}{\input{fix-poset-ex.pstex_t}}
\caption{The posets $[\s, \t]$ and $[\lsr, \t]$, with $\s = 321$, $\t =
634521$, and $\lsr = \la642\ra$.  In the
latter poset, the letters in the occurrence $\lsr$ are marked.}\label{fig:posets}
\end{figure}

When no confusion will arise, the M\"obius function of the interval $[\lsr, \t]$
will be denoted
$$\mu(\lsr,\t) := \mu_{[\langle \sigma \rangle, \tau]}(\lsr, \t).$$

Definition~\ref{defn:occurrence interval} will be helpful for
  analyzing the M\"obius function of a generic interval $[\s,\t]$ in
  $\p$ because the M\"obius function can be computed as an alternating
  sum of chains according to length.  In any chain from $\t$ to $\s$
  in $\p$, the final element must be a $\s$-pattern.  Thus this chain
  highlights a particular occurrence of $\s$ in $\t$ and removes the
  other letters successively.  Of course, we must exercise caution
  because there may be some overcounting of these chains.  One instance in
  which there can be no overcounting is when $\s(\t) = 1$, and hence
  $[\s,\t] = [\lsr, \t]$.

\section{Pairs with interval blocks}\label{section:interval blocks}

Throughout this section, fix $\s \in \sym_k$ and $\t \in \sym_n$ such
that $\s \le \t$ in $\p$, and fix an occurrence $\lsr$ of $\s$ in
$\t$.

\begin{thm}\label{thm:interval blocks - occurrence}
  Let $\s \in \sym_k$ and $\t \in \sym_n$, with $\s\le\t$.  If the
  pair $(\lsr, \t)$ has an interval block, then
$$\mu(\lsr,\t) = 0.$$
\end{thm}

\begin{proof}
  We prove the theorem by induction on $n-k$, the rank of the interval
  $\lsrt$.  Since $(\lsr, \t)$ has an interval block, $n-k$ must be at
  least 2. In the case $n-k=2$, which serves as the basis of the
  induction, the pair $(\lsr, \t)$ has precisely one interval block (with
  two elements), which we call $I$.  This implies that the poset
  $\lsrt$ is a three element chain, where the middle element is
  obtained by removing one of the two letters of $I$ from $\t$.  The
  M\"obius function of such a chain is 0, as desired.

  Assume that the result holds for all
  intervals of rank less than $n-k$, where $n-k\ge3$.  Let $I$ be an
  interval block in the pair $(\lsr,\t)$, where $[\lsr,\t]$ has rank
  $n-k$.  Let $\pi$ be the permutation obtained from $\t$ by deleting
  all but one letter from $I$, but otherwise leaving $\t$ intact.
  For each permutation $\rho\in\lsrt$ such that
  $\rho\not\le\pi$, the pair $(\lsr,\rho)$ must contain an interval
  block because $\rho$ must contain at least two letters originating
  from $I$.

Because $[\lsr,\pi]$ is a
  closed interval in $[\lsr,\t]$, the sum of the values of the
  M\"obius function over this interval is zero.  Thus
\begin{eqnarray}
\nonumber \mu(\lsr,\t) &=& - \sum_{\rho \in [\langle \sigma \rangle, \tau)} \mu(\lsr, \rho)\\
\nonumber &=& - \sum_{\rho \in [\langle \sigma \rangle, \pi]} \mu(\lsr, \rho) - \sum_{\substack{\rho \in [\langle \sigma \rangle, \tau) \\ \rho \not\le \pi}} \mu(\lsr, \rho)\\
\label{eqn:interval block induction}&=& - \sum_{\substack{\rho \in [\langle \sigma \rangle, \tau) \\ \rho \not\le \pi}} \mu(\lsr, \rho).
\end{eqnarray}
Now, each $\rho$ in the summation in Equation~\eqref{eqn:interval
  block induction} is strictly smaller than $\t$, so the rank of
$[\lsr,\rho]$ for each such $\rho$ is strictly less than $n-k$.
Moreover, the pair $(\lsr,\rho)$ contains an interval block.  Hence,
by the inductive hypothesis, $\mu(\lsr,\rho)=0$ for all such $\rho$,
which implies that $\mu(\lsr,\t)=0$, completing the proof.
\end{proof}

As discussed earlier, the poset $[\lsr, \t]$ will be a helpful tool
for the analysis of intervals in $\p$.  In fact,
Theorem~\ref{thm:interval blocks - occurrence} can be translated
readily into a statement about intervals in $\p$, as shown in the
following theorem.  The proof of this result is entirely analogous to
that of the previous theorem.

\begin{thm}\label{thm:interval blocks}
If the pair $(\s, \t)$ has an interval block, then $$\mu(\s,\t) = 0.$$
\end{thm}

\begin{proof}
Consider $\s \in \mathcal{S}_k$ and $\t \in \mathcal{S}_n$, with $\s \le \t$ in $\p$.

We prove the theorem by induction on $n-k$, the rank of the interval $[\s, \t]$. Since $(\s, \t)$ has an interval block, $n-k$ must be at least $2$.  In the case $n-k = 2$, which serves as the basis of the induction, the pair $(\s, \t)$ as precisely one interval block (with two elements), which we call $I$.  This implies that the poset $[\s, \t]$ is a three element chain, where the middle element is obtained by removing one of the two letters of $I$ from $\t$.  The M\"obius function of such a chain is $0$, as desired.

Assume that the result holds for all intervals of rank less than $n-k$, where $n-k \ge 3$.  Let $I$ be an interval block in the pair $(\s, \t)$, where $[\s, \t]$ has rank $n-k$.  Let $\pi$ be the permutation obtained from $\t$ by deleting all but one letter from $I$, but otherwise leaving $\t$ intact.  For each permutation $\rho \in [\s, \t]$ such that $\rho \not\le \pi$, the pair $(\s, \rho)$ must contain an interval block because $\rho$ must contain at least two letters originating from $I$.

Because $[\s,\pi]$ is a closed interval in $[\s,\t]$, the sum of the values of the M\"obius function over this interval is zero.  Thus
\begin{eqnarray}
\nonumber \mu(\s,\t) &=& - \sum_{\rho \in [\sigma, \tau)} \mu(\s, \rho)\\
\nonumber &=& - \sum_{\rho \in [\sigma, \pi]} \mu(\s, \rho) - \sum_{\substack{\rho \in [\sigma, \tau) \\ \rho \not\le \pi}} \mu(\s, \rho)\\
\label{eqn:interval block induction2}&=& - \sum_{\substack{\rho \in [\sigma, \tau) \\ \rho \not\le \pi}} \mu(\s, \rho).
\end{eqnarray}
Now, each $\rho$ in the summation in Equation~\eqref{eqn:interval
  block induction2} is strictly smaller than $\t$, so the rank of
$[\s,\rho]$ for each such $\rho$ is strictly less than $n-k$.
Moreover, the pair $(\s,\rho)$ contains an interval block.  Hence,
by the inductive hypothesis, $\mu(\s,\rho)=0$ for all such $\rho$,
which implies that $\mu(\s,\t)=0$, completing the proof.
\end{proof}

There are several things to note about Theorem~\ref{thm:interval
  blocks}.  First, the result does not hold if there is an interval
block in $\t$ that is disjoint from some, but not all, occurrences of
$\s$ in $\t$.  For example, let $\s = 12$ and $\t = 3412$.  The
interval block $12$ is disjoint from the occurrence $\lsr_1 = \la34\ra$, but
not from the occurrence $\lsr_2 = \la12\ra$.  However, if one considers the
poset $[\lsr_1, \t]$, then $\mu_{[\langle \sigma \rangle_1,
  \tau]}(\lsr_1, \t) = 0$, whereas $\mu_{[\sigma,
  \tau]}(\s,\t)=1$.  

Also, the converse of Theorem~\ref{thm:interval blocks} is false.  For
example, let $\s = 1$ and $\t = 123$.  The pair $(\s,\t)$ does not
have an interval block, but $\mu(\s,\t) = 0$.  Thus, the property
$\mu(\s,\t) = 0$ cannot be characterized completely by the presence of
an interval block.

\section{Intervals describing a fixed occurrence}\label{section:fixed occurrence}
  
Throughout this section, fix $\s \in \sym_k$ and $\t \in \sym_n$ such
that $\s \le \t$ in $\p$, and fix an occurrence $\lsr$ of $\s$ in
$\t$.

In this section we give two results showing that the M\"obius
function of a class of intervals is $1$ or $-1$, the sign depending on
only on the rank of the interval.  The first of these,
Theorem~\ref{thm:boolean algebra}, deals with intervals that are
boolean algebras.  That, in turn, is the basis for an extension to a more general result in Corollary~\ref{cor:single occurrence}.

\begin{defn}
Let $\lsr$ be an occurrence of $\s$ in $\t$.  A \emph{region} in $\t$ (with respect to $\s$) is a maximal consecutive substring of the complement of $\lsr$ in $\t$.  
\end{defn}

We now define a quality of a pair $(\lsr, \t)$, which informally means that every pair of distinct letters in $\t \setminus \s$ is separated either in value or in position by some letter of $\lsr$.  This quality is defined more precisely as follows.

\begin{defn}\label{defn:separated}
The pair $(\lsr, \t)$ is \emph{separated} if for all $x < y$ that are in the same region of $\t \setminus \lsr$, there exists a $j$ such that $x < \la \s(j) \ra < y$.
\end{defn}

\begin{ex}
  The pair $(\lsr, \t) = (\langle 653\rangle, 146253)$ is separated: although $1$ and $4$ lie in the same region of $\t \setminus \lsr$, we have $1 < \la \s(1) \ra = \la3\ra < 4$.

  The pair $(\lsr, \t) = (\langle 764 \rangle, 1357264)$ is not
  separated because $1$ and $3$ belong to the same region of $\t
  \setminus \lsr$, but the only letter between them in value, the 2,
  does not lie in $\lsr$.
\end{ex}

\begin{thm}\label{thm:boolean algebra}
The interval $[\lsr,\t]$ is
  boolean if and only if the pair $(\lsr,\t)$ is separated.  In this
  case, therefore, $\mu(\lsr,\t)=(-1)^{n-k}$, where $k$ and $n$ are the
  ranks of $\s$ and $\t$, respectively.
\end{thm}

\begin{proof}
  The only letters that can be deleted from $\t$ in the interval
  $[\lsr,\t]$ are elements of $\t \setminus \lsr$, since letters of
  the occurrence $\lsr$ must be intact in each element of the
  interval.  Suppose that $\s \in \sym_k$ and $\t \in \sym_n$, so $[\s,\t]$ has
  rank $n-k$.  It is boolean if and only if the letters of $\t
  \setminus \lsr$ can be removed in any order, always yielding
  distinct permutations.  This happens if and only if at no point are
  there letters $x$ and $x+1$ adjacent in $\rho \setminus \lsr$ for
  some $\rho \in [\lsr,\t]$.  This is equivalent to the pair
  $(\lsr,\t)$ being separated.
\end{proof}
It is well known that the M\"obius function of a boolean algebra is
$\pm1$, the sign depending only on rank.  This yields the following
corollary.

\begin{cor}
  If $\s(\t) = 1$ and the pair $(\lsr, \t)$ is separated, then
  $\mu(\s,\t) = (-1)^{n-k}$, where $k$ and~$n$ are the ranks of $\s$
  and $\t$, respectively.
\end{cor}

It should be noted that $[\s,\t]$ can be boolean even if $\sigma(\tau)
> 1$.  For example, the intervals $[123,1324]$ and $[12,3412]$ are
both boolean, even though $123(1324) = 12(3412) = 2$.  It would thus
be interesting to characterize boolean intervals.

In Section~\ref{section:interval blocks}, pairs $(\lsr, \t)$ with
interval blocks were analyzed and shown to satisfy $\mu(\lsr, \t) =
0$.  Moreover, an interval $[\s,\t] \subset \p$ where the pair
$(\s,\t)$ has an interval block satisfies $\mu(\s,\t) = 0$.  Now we
examine situations where $\s \le \t$ in $\p$ and there is an
occurrence $\lsr$ such that the pair $(\lsr, \t)$ is interval free.
Under certain conditions, the M\"obius function for such a poset has
the nice form $\mu(\lsr, \t)=(-1)^r$, where $r$ is the rank of
$\lsrt$.  We believe that this property holds for a large class of
intervals and give one such subclass here.

A few preliminaries are necessary before proving Theorem~\ref{thm:fixed occurrence}.

\begin{defn}\label{defn:lsrtt}
The subposet $\ifset{\lpr}{\rho}$ of $[\lpr, \rho]$ consists of all $\omega \in
[\lpr,\rho]$ such that the pair $(\lpr,\omega)$ is interval free.  In
particular, it is always the case that $\lpr \in \ifset{\lpr}{\rho}$, and $\rho \in \ifset{\lpr}{\rho}$ if the pair $(\lpr, \rho)$ is interval free.
\end{defn}

\begin{rem}\label{rem:removing 0s}
  By Theorem~\ref{thm:interval blocks - occurrence}, the poset
  $\lsrtt$ is obtained by removing from $[\lsr, \t]$ only elements
  $\pi$ with $\mu(\lsr, \pi) = 0$.  When the pair $(\lsr, \t)$ is interval free, we are thus
  replacing $[\lsr, \t]$ by a subposet, $\lsrtt$, which has the
  same M\"obius function as $[\lsr, \t]$.
\end{rem}

\begin{defn}\label{def:rp}
  A subset $S$ of a graded poset $P$ has the \emph{rank property} if
  it contains as many elements of even rank in $P$ as it does of odd
  rank.  Such an $S$ is called \emph{\rp}.  If $S$ is all of $P$ we
  use the same terminology for $P$ itself.
\end{defn}

The following lemma suggests the relevance of the rank property.

\begin{lem}\label{lem:RP}
An \rp\ poset $P$ has the following properties.
\begin{enumerate}
\item If $X \subseteq P$ is \rp\ as well, then $P\setminus X$ is \rp.
\item\label{lemRPb} If $P$ is bounded, with minimal element $\hat{0}$ and maximal
  element~$\hat{1}$, and if $\mu_P(\hat{0},x) = (-1)^{\rk(x)}$ for all
  $x \in [\hat{0}, \hat{1})$, then $\mu_P(\hat{0}, \hat{1}) =
  (-1)^{\rk(P)}$.
\end{enumerate}
\end{lem}

\begin{proof}
Part (a) is obvious.

By Definition~\ref{defn:mu defn},
    \begin{equation}\label{eqn:rank symmetry}
\mu(\hat{0}, \hat{1}) = -r_0 + r_1 -r_2 + \cdots - (-1)^{\rk(P)-1}r_{\rk(P)-1},
\end{equation}
where $r_i$ is the number of elements of rank $i$ in $P$.  The poset $P$ is {\rp}, so
$$r_0 + r_2 + r_4 + \cdots = r_1 + r_3 + r_5 + \cdots,$$
with one sum going to $r_{\rk(P)-1}$ and the other to $r_{\rk(P)}$.
This is equivalent to
$$ (-1)^{\rk(P)}r_{\rk(P)}= -r_0 + r_1 - r_2 + \cdots - (-1)^{\rk(P)-1}r_{\rk(P)-1}.$$
Thus equation~\eqref{eqn:rank symmetry} simplifies to $\mu(\hat{0},
\hat{1}) = (-1)^{\rk(P)}r_{\rk(P)}$.  Since $r_{\rk(P)} = 1$, this
means that $\mu(\hat{0}, \hat{1}) = (-1)^{\rk(P)}$, proving part (b).
\end{proof}

We will show that under particular conditions, the poset
$\ifset{\lsr}{\t}$ is \rp, regardless of whether the pair $(\lsr, \t)$
is interval free.  When the pair is interval free, and thus $\t \in
\lsrtt$, this will imply that $\mu(\lsr, \t) = (-1)^{n-k}$ in both
$[\lsr, \t]$ and $\ifset{\lsr}{\t}$.

Recall Definition~\ref{defn:occurrence interval} of the set $\C = \t
\setminus \lsr$: $C$ consists of those values of $[1,n]$ that are not
part of the occurrence $\lsr$.  This set holds the key to the proof
because we have fixed the occurrence $\lsr$ of $\s$ in $\t$, so
elements of $\lsrt$ differ only in the letters of $\C$.

Any permutation $\rho\in\lsrt$ can be identified with the subset of
letters $\rho \setminus \lsr \subseteq \C$.  However, there may be
distinct subsets $S$ and $S'$ of $\C$ for which the permutations $\lsr
+ S$ and $\lsr + S'$ are order isomorphic, having the occurrence
$\lsr$ in the same places.  Of course, two such permutations
correspond to the same element of $\lsrt$.  When discussing something
of the form $\lsr + S$, where the letters of $S$ are added to the
occurrence $\lsr$, it will be understood that this indicates the order
isomorphic element in $\lsrt$ (or in $\lsrtt$ as appropriate), having
the occurrence $\lsr$ in the same positions.

We now come to the definitions of the conditions under which the
M\"obius function of an interval free pair $(\s,\t)$ can be shown to
be either $1$ or $-1$, depending only on the rank of the interval.
Note that there is a relationship between the following definition and
the notion of ``separated'' (see Definition \ref{defn:separated}).

\begin{defn}
  Two letters $a$ and $b$, with $a<b$,
  in the complement of $\s$ in $\t$ are \emph{similar} if they are in
  the same region and there is no $c\in\lsr$ such that $a<c<b$.  A
  \emph{group of similar letters} in the complement consists of two or
  more similar letters.
\end{defn}

As an example, if $\t=357128469$ and $\lsr = \langle184\rangle$, then the regions are $357-2-69$.  The letters 3 and 5
are not similar, since $4 \in \lsr$.  The letters 5 and 7 are
similar, because the only letter between them in value, namely 6,
belongs to the complement of $\lsr$.

\begin{thm}\label{thm:fixed occurrence}
  Fix $\s \in S_k$ and $\t \in S_n$, where $k<n$, and let $\lsr$ be
  an occurrence of $\s$ in $\t$.  If there is at most one maximal
  group of similar letters in all the regions of $\t$, then
  $\ifset{\lsr}{\t}$ is \rp.
\end{thm}

\begin{proof}
  Suppose that there is no maximal group; that is, no two letters in
  the same region are similar.  Then the pair $(\s,\t)$ is separated,
  so, by Theorem~\ref{thm:boolean algebra}, the interval $\lsrtt$ is a
  nontrivial boolean algebra, which is \rp.

Now suppose that there is a unique maximal group of similar letters,
and that it consists of $m \ge 2$ letters.  We prove the result by
induction on $(n-k)-m$.

If $n-k = m$, then the entire complement $\t \setminus \lsr$ is this
group of $m$ similar letters, which forms an interval block.  Thus
$\ifset{\lsr}{\t}$ consists of two elements: the element $\lsr$ and
the element obtained by adding to $\lsr$ a single letter from the
complement $\t \setminus \lsr$.  This poset is obviously \rp.

Now assume the result for all $(n-k)-m < r$, where $r \ge 1$ and
suppose that $(n-k)-m = r$.  This means that in addition to the group
of $m$ elements, there are $r$ other letters in $\t \setminus \lsr$,
none of which is similar to any other letter in its respective
region. Suppose that $x$ is one of these letters.  Then the poset
$\ifset{\lsr}{\t}$ is isomorphic to the disjoint union of the posets
$\ifset{\lsr}{\t \setminus \{x\}} \sqcup \ifset{\lsr+x}{\t}$, based on
whether an element of $\ifset{\lsr}{\t}$ includes the letter $x$ or
not.  Each of these intervals satisfies the hypotheses of the theorem,
and by induction we know that they are both \rp.  Therefore
$\ifset{\lsr}{\t}$ is \rp\ as well.
\end{proof}

\begin{cor}\label{cor:single occurrence}
  If $\s$ and $\t$ satisfy the hypotheses of
  Theorem~\ref{thm:fixed occurrence}, and the pair $(\s,\t)$ is interval
  free, then $\mu(\s,\t)=(-1)^{n-k}$.
\end{cor}

\begin{proof}
This follows from Lemma~\ref{lem:RP} and Theorem~\ref{thm:fixed occurrence}.
\end{proof}

The hypotheses of Theorem \ref{thm:fixed occurrence} cannot be
weakened to encompass all interval free pairs $(\s,\t)$.  The smallest
$\t$ for which there are counterexamples $\lsrt$ with
$\mu(\s,\t)\ne\pm1$ have length 10, where the shortest $\s$ have
length 3.  The only counterexample in that case, up to trivial
symmetries, is when $\s=321$ and $\t=2,5,1,7,3,10,4,6,9,8$ (we use
commas here to clarify which values are two-digit numbers).  The only
occurrence of $\s$ in $\t$ here is $\langle 10,9,8 \rangle$, so the
complement is $25173-46$, which has no interval blocks, but
$\mu(\s,\t)=0$.

Another counterexample, where the M\"obius function has value $2$, is
$$
\s=2341,~~~ \t=2,3,8,1,6,12,4,10,5,9,7,11,2.
$$  
Here the only occurrence of $2341$ is $\langle2381\rangle$.  Thus the complement
is a single region of eight letters.

\section{Open problems}\label{section:multiple occurrences}

We have gathered data to support the following two conjectures.

\begin{conj}
  Suppose that $\t$ avoids the pattern $132$ (or, equivalently, $312$,
  or $213$, or $231$).  Then $\mu(1,\t)$ is either $0$, $1$, or $-1$.
\end{conj}
\begin{conj}\label{conj:132}
  Suppose that $\t$, and hence the entire interval $[\s,\t]$ for any
  $\s$, avoids the pattern $132$ (or, equivalently, $312$, or $213$,
  or $231$).  Then we have $|\mu(\s,\t)| \le \s(\t)$.
\end{conj}

Because of the symmetries between the patterns $132$, $231$, $213$,
$312$, if Conjecture \ref{conj:132} is true for any one of them, it is
true for each of them.  Thus, the combined contrapositives of the
corresponding four conjectures say that if $|\mu(\s,\t)|$ exceeds
$\s(\t)$, then $\t$ must contain all of the patterns $132$, $231$,
$213$ and $312$.

Observe that if Conjecture \ref{conj:132} holds, then, for intervals
$[\s,\t]$ where $\s$ occurs precisely once in $\t$, the M\"obius
function must be $0$, $1$ or $-1$.

Additionally, recalling Section~\ref{section:fixed occurrence}, it would be interesting to characterize the boolean intervals in the poset of permutation patterns.

\section*{Acknowledgment}

We are indebted to Eric Babson for very helpful discussions, and to an
anonymous referee for many valuable suggestions.

\end{document}